\documentclass[11pt,a4paper]{article}

\usepackage[a4paper,body={15.6cm,23cm}]{geometry}

\usepackage[intlimits]{amsmath}
\usepackage{amsfonts}
\usepackage{mathrsfs}

\usepackage{bbm}
\usepackage{dsfont}

\usepackage{theorem}

\usepackage{QED}

\newcommand{\N}{\mathbbm{N}}	
\newcommand{\Z}{\mathbbm{Z}}	
\newcommand{\R}{\mathbbm{R}}	

\theoremstyle{plain}
\theoremheaderfont{\bf}
\newtheorem{definition}{Definition}[section]
\newtheorem{proposition}[definition]{Proposition}
\newtheorem{lemma}[definition]{Lemma}

\newtheorem{theorem}[definition]{Theorem}

\begin{document}

\title{Entropic repulsion for a class \\of Gaussian interface models in high dimensions}
\author{Noemi Kurt\thanks{E-mail: noemi.kurt@math.unizh.ch, Phone: ++41 44 635 58 43, Fax: ++41 44 635 57 05.}\\ Universit\"at Z\"urich,
Institut f\"ur Mathematik,\\ Winterthurerstr. 190, CH-8057 Z\"urich} 
\date{February 6, 2006}
\maketitle

\begin{abstract}
Consider the centered Gaussian field on the lattice $\mathbb{Z}^d,$ $d$ large enough, with covariances given by the inverse of $\sum_{j=k}^K q_j(-\Delta)^j,$ where
$\Delta$ is the discrete Laplacian and $q_j \in \mathbb{R},k\leq j\leq K,$ the $q_j$ satisfying certain additional conditions. We extend a
previously known result to show that the probability that all spins are nonnegative
on a box of side-length $N$ has an exponential decay at rate of order $N^{d-2k}\log{N}.$ The constant is given in terms of a higher-order capacity of the unit cube, analogous
to the known case of the lattice free
field. This result then allows us to show that, if we condition the field to stay positive in the $N-$box, the local sample mean of the field is pushed to a height of
order $\sqrt{\log N}.$ \\

\noindent
Keywords: Random Interfaces, Entropic Repulsion, Gaussian Fields.
\end{abstract}

\section{Introduction and results}
We study the entropic repulsion of a class of real valued Gaussian random fields $%
\varphi =\left\{ \varphi _{x}\right\} _{x\in \mathbbm{Z}^{d}},$ which can be interpreted as a $d$-dimensional (discrete) interface
in a $\left( d+1\right) $-dimensional space. Entropic repulsion refers to
the fact that the presence of a wall forces the random surface to move
away from the wall, in order to gain space for local fluctuations (cf. \cite{LebowitzMaes}). In our case, the wall is simply the $d$-dimensional
coordinate hyperplane, and the effect of the wall is described by requiring
the field $\left\{ \varphi _{x}\right\} $ to be positive in a certain region.

A basic object to study is the asymptotics of the probability $P(
\varphi _{x}\geq 0,\ x\in V) $, $V\subset \mathbbm{Z}^{d}$ finite, when
$V\uparrow \mathbbm{Z}^{d}.$ Its behaviour is well understood in the
case where $\left\{ \varphi _{x}\right\} $ is the (lattice) Gaussian free
field in dimension $d\geq 3$ (see \cite{BDZ}), which is the Gibbs measure with formal
Hamiltonian $H( \varphi ) =\sum_{\left\vert x-y\right\vert
=1}( \varphi _{x}-\varphi _{y}) ^{2}.$ The free field has a
simple random walk representation of the covariances, which enables one to
calculate various conditional distributions in an easy way. The main aim of
the present paper is to extend the analysis to a class of fields with a more
general Hamiltonian which includes the so-called \textquotedblleft membrane
models\textquotedblright . The crucial difference is that these models do not
possess a random walk representation.

The only mathematically rigorous result on these models we are aware of is
the paper by Sakagawa \cite{Sakagawa}, who derived lower and upper bounds for $P\left(
\varphi _{x}\geq 0,\ x\in V\right) .$ His bounds however don't match. We
derive here an upper bound which asymptotically matches Sakagawa's lower
bound, and therefore we prove that his lower bound gives the correct leading
order of the asymptotics. This first result then enables us to compute the exact height to which the average of the field is
pushed by the wall.

Let $\Omega=\R^{\Z^d}.$ We consider the (formal) Hamiltonian $H:\Omega\to\R$ given by 
\begin{eqnarray}\label{Hamilton} H(\varphi)=\sum_{j=1}^Kq_j\sum_{x\in \Z^d}((-\Delta)^{\frac{j}{2}}\varphi_x)^2\,,\end{eqnarray}
with $q_j \in \R,$ for $1\leq j \leq K.$ The discrete Laplacian $\Delta$ is the operator on $L^2(\Z^d)$ defined by $\Delta
f(x)=\left(\frac{1}{2d}\sum_{y\in \Z^d, |x-y|=1}f(y)\right)-f(x).$ If $j$ is odd, we set
$\sum_{x\in \Z^d}((-\Delta)^{\frac{j}{2}}\varphi_x)^2=\sum_{i=1}^d\sum_{x\in \Z^d}((-\Delta)^{\frac{j-1}{2}}\nabla_i\varphi_x)^2,$ where $\nabla_i$ denotes the discrete
gradient in the $i-$th direction. The
free field is thus the special case $K=1,\, q_1=1.$ 
We make the following
assumptions:
\begin{itemize}
\item[(a)] $d\geq 2k+1,$ where $k=\min\{j: q_j\neq 0\}\,,$
\item[(b)] $q=\{q_j\}_{1\leq j\leq K} \in \R^K$ satisfies $q(r):=\sum_{j=k}^Kq_jr^j>0$ for $0<r\leq 2\,.$
\end{itemize}
Under assumptions $(a)$ and $(b),$ the infinite-volume Gibbs
measure corresponding to $H$ exists (see \cite{Sakagawa}, Section 2). It can be described as follows: 
For $\varepsilon \geq 0$ and $x, y\in \Z^d$ set $J_\varepsilon(x,y)=q(\varepsilon I -\Delta)(x,y),$ where $I$ is the identity matrix on $\Z^d,$
and $\Delta$ the matrix Laplacian defined by
$$\Delta(x,y)=\left\{\begin{array}{rl} -1 & \text{if}\, x=y\,,\\
\frac{1}{2d} &  \text{if}\, |x-y|=1\,,\\
0&  \text{otherwise .}\end{array}\right.$$  The
above assumptions ensure that, for $\varepsilon$ small enough, the matrix $J_\varepsilon$ is positive definite with positive definite inverse
$J_\varepsilon^{-1}.$ We set $J(x,y)=J_0(x,y)$ and $G(x,y)=J^{-1}(x,y).$ From \cite{Georgii}, Chapter 13, we know that the centred Gaussian field with covariance matrix $G$
exists. We denote its law by $P.$ It is characterised by the following DLR-equation as an infinite-volume
Gibbs measure, and corresponds to the Hamiltonian (\ref{Hamilton}):
\begin{eqnarray}\label{DLR}P(\:\cdot\:|{\cal{F}}_{\{x\}^c})(\varphi)={\cal{N}}(-J(0,0)\cdot\sum_{y\neq x}J(x,y)\varphi_y,\;J(0,0)^{-1}) \quad
P-a.s.\,,\end{eqnarray}
where we use the notation ${\cal{F}}_A=\sigma (\varphi_y: y\in A)$ for the $\sigma-$field generated by $\{\varphi_y:y\in A\},$ $A\subset
\Z^d,$ and ${\cal N}(\mu, \sigma^2)$ is the Gaussian distribution with mean $\mu$ and variance $\sigma^2.$ If $G(x,y)\geq 0$ for all $x,y,$ then $P$ satisfies the
FKG-inequalities (see \cite{Pitt}). We make the additional assumption
 \begin{itemize}
 \item[(c)] There exists a sequence $\{\varepsilon_n\}_{n\in \N}$ of positive numbers such that
 $\lim_{n\to\infty}\varepsilon_n=0$ and
 $J_{\varepsilon_n}^{-1}(x,y)\geq 0$ for all $n\in \N$ and all $ \;x,y\in \Z^d.$
 \end{itemize}
Throughout the paper we will always assume that $(a), (b)$ and $(c)$ hold. 

Now set $V=[-1,1]^d$ and $ V_N=NV\cap \Z^d.$ We consider the entropic repulsion event
$$\Omega_N^+=\left\{\varphi\in \Omega: \varphi_x\geq 0 \; \forall x\in V_N\right\}.$$
 Theorem 2.1 of \cite{Sakagawa} states that there exist constants $C_1, C_2>0$ such that
\begin{eqnarray}\label{thmSaka}
-C_1\leq\liminf_{N\to\infty}\frac{1}{N^{d-2k}\log{N}}\log{P(\Omega_N^+)}\leq
\limsup_{N\to\infty}\frac{1}{N^{d-2k}\log{N}}\log{P(\Omega_N^+)}\leq -C_2
\end{eqnarray}
holds. Moreover, the constant $C_1$ has been identified to be $C_1=2k\cdot q_k\cdot G(0,0)\cdot C_k(V),$ where
\begin{eqnarray}\label{constSaka} C_k(V)=\inf\left\{\frac{1}{(2d)^k}\int_{\R^d}|(-\nabla)^kh|^2 dx;\,\, h\in H^k(V), \, h\geq 1 \;\text{on}\;
V\right\}\end{eqnarray} is the $k-$th order capacity of the unit cube $V,$ $k$ being the minimal degree of the polynomial $q.$ This lower bound was proved using a relative entropy argument
and the FKG-property of $P$. Assumption (c) above is necessary for this proof. In the case of the
free field, the  statement \ref{thmSaka} was proved before in \cite{BDZ}. There it was shown in addition that in this case the constants $C_1$ and $C_2$ of the upper and the
lower bound coincide. Our first
result shows that this is still
true for our model:
\begin{theorem} \label{thmer}

For $d\geq 2k+1,$ 

\begin{eqnarray}\label{er}\lim_{N\to\infty}\frac{1}{N^{d-2k}\log N}\log P\left(\Omega_N^+\right)=-2k\cdot q_k\cdot G\cdot C_k(V)\,,\end{eqnarray} 

where $C_k(V)$ is given by \rm(\ref{constSaka}), and $G=G(0,0).$ 
\end{theorem}

In the next section, we will prove the upper bound of (\ref{er}). Together with (\ref{thmSaka}) and (\ref{constSaka}) this proves Theorem \ref{thmer}. Thus the decay of
$P(\Omega_N^+)$ for $k\geq 2$ is completely analogous to the case $k=1.$ Using Theorem \ref{thmer}, we can then prove
the height estimate for the averaged field:

\begin{theorem}\label{thmheight}
Let $\varepsilon>0$ and $\eta>0.$ Then
\begin{eqnarray}\lim_{N\to\infty}\sup_{z\in V_N,\atop V_{N,\varepsilon}(z)\subset V_N}P\left(\left|\frac{\overline{\varphi}_{N,\varepsilon}(z)}{
\sqrt{\log
N}}-\sqrt{4kG}\right|\geq \eta\;\Bigg| \,\Omega_N^+\right)=0\,,\label{height}\end{eqnarray}
where $\overline{\varphi}_{N,\varepsilon}(z)=\frac{1}{|V_{N,\varepsilon}(z)|}\sum_{x\in V_{N,\varepsilon}(z)}\varphi_x$ and
$V_{N,\varepsilon}(z)=\{x\in V_N: \max_{1\leq i\leq d}|x_i-z_i|\leq \varepsilon N\}.$
\end{theorem}

The lower bound for this height estimate was obtained in \cite{Sakagawa}. Our exact result in Theorem \ref{thmer} allows us now to give the
correct upper bound. This means that, as expected, the local sample mean of the field is pushed to $\sqrt{4kG\cdot \log{N}}$ by the hard wall.

\section{Proof of the upper bound in Theorem \ref{thmer}}

We follow a strategy introduced in \cite{BertGiac02}, which was used in \cite{Giacomin} for the case $k=K=1, \, q_k=1.$ The idea is to use a conditioning argument on larger
boxes than those of the proof of \cite{Sakagawa}. The main difficulty -- when trying to follow the proof for the free field -- arises when considering the
expectations of $\varphi_x$ conditioned on the boundary of a box of side-length $L.$ While in the case of the harmonic crystal, we
know by the random walk representation, that on $\Omega_N^+$ the conditional expectations are nonnegative, in our more general case they
can be
strictly negative. We overcome this difficulty by estimating the proportion of conditional expectations that are of order $-N^\lambda,\,
\lambda \in \N.$ Then we prove that this proportion is negligible if we let $N$ tend to infinity.\\

\it Proof of Theorem \ref{thmer}, the upper bound. \rm 
Fix a natural number $L>K+1$ such that $L-K$ is even, and let $\overline{\Lambda}=(L, L,...,L)+L\Z^d.$ For $x\in \overline{\Lambda}$
denote by $\partial B(x)=\{y\in \Z^d: \max_{i=1,...,d}|x_i-y_i|\in [\frac{L-K}{2},\frac{L+K}{2}]\}$ the boundary of the box $B(x):=\{y\in
\Z^d: \max_{i=1,...,d}|x_i-y_i|<\frac{L-K}{2}\}.$ Let $\tilde{\Lambda}=\{x\in \overline{\Lambda}:\partial B(x)\subset V_N\}$ and
$\Lambda=\cup_{x\in \tilde{\Lambda}}\partial B(x).$\\

Since $J(x,y)=0$ for $|x-y|>K,$ the field $\{\varphi_x\}_{x\in \tilde{\Lambda}}$ is Markovian, in the sense that $P(\; \cdot \; |\,{\cal
F}_{B(x)^c})=P(\;\cdot\;|\,{\cal F}_{\partial B(x)})$ for all $x\in \tilde{\Lambda},$ and thus (see \cite{Georgii}, Proposition 13.13), under
$P(\;\cdot\;|\,{\cal F}_{B(x)^c}),$ the $\varphi_x,\,x\in \tilde{\Lambda},$ are independent normally distributed random variables. For the
mean and the variance we write 
$$m_x=E\left(\varphi_x\,|\,{\cal F}_{B(x)^c}\right)\quad \text{and} \quad G_L={\rm var}\left(\varphi_x\,|\,{\cal F}_{B(x)^c}\right)$$
respectively. Note that $\lim_{L\to\infty}G_L=G$ (see \cite{Georgii}, Section 13.1).
For any subset $A$ of $\Z^d$ let $\Omega_A^+$ denote the event $\{\varphi_x\geq 0 \quad \forall x\in A\}.$ Because of the independence we have 
\begin{eqnarray}\label{unabh}P\left(\Omega_N^+\right)\leq P\left(\Omega^+_{\Lambda}\cap\Omega^+_{\tilde{\Lambda}}\right)\leq E\left[\prod_{x\in \tilde{\Lambda}}P(\varphi_x\geq 0\,|\,{\cal F}_{\partial
B(x)})\cdot 1_{\Omega_\Lambda^+}\right].\end{eqnarray}
As in \cite{Giacomin}, we use a decomposition of $V$ on a larger scale: Let $\theta >0,$ $r\in \R^d$ and set $A_r=r+[0,\theta)^d,$ and
$I=\{r\in \theta \Z^d : A_r\subset V, \, \partial A_r\cap \partial V=\emptyset\}.$ Set $\tilde{B}_r=NA_r\cap \tilde{\Lambda},$ the
box containing the centres of the smaller boxes $B(x),$ with $x\in NA_r.$ Note that
$B:=|\tilde{B}_r|=O(N^d).$ \\

Let $0<\delta<1$ and $0<\gamma<1.$ For $\kappa>0,$ define $a_N=\sqrt{4k(G-\kappa)\log{N}}$ and consider the following events:
\begin{eqnarray*} E_{\delta, \kappa}&=&\left\{\varphi: \text{there is} \,r\in I \,\text{such that} \:\big|\{x\in \tilde{B}_r: m_x\leq
a_N\}\big|\geq \delta B\right\},\\
E_{\delta}^{-\lambda}&=&\left\{\varphi: \text{there is}\, r \in I \, \text{such that } \: \big|\{x\in \tilde{B}_r: m_x\leq
-N^{\lambda} \}\big|\geq  \frac{\delta}{N^{2k+(2-\gamma)\lambda-\gamma}}B\right\}, 
\end{eqnarray*}
 where $\lambda$ is a nonnegative integer. Note that for $\lambda\geq\lambda_{\max}=\left\lfloor(d-2k+\gamma)/(2-\gamma)\right\rfloor+1$ (where
 $\big\lfloor\cdot\big\rfloor$ denotes the integer part), we have
 $N^{-2k-(2-\gamma)\lambda+\gamma}\delta B<1.$ For these $\lambda$s, $E_{\delta}^{-\lambda}$ is $\{\varphi: \text{there is}\, r \in I \, \text{such that } \, \{x\in \tilde{B}_r: m_x\leq
-N^{\lambda} \}\neq\emptyset\},$ and these $E^{-\lambda}_\delta$s are all contained in $E^{-\lambda_{\max}}_\delta.$ Set
$$F_{\delta}=\bigcup_{\lambda=0}^{\lambda_{\max}}E_{\delta}^{-\lambda}.$$

The estimate (\ref{unabh}) now gives
$$P(\Omega_N^+)\leq E\left[\prod_{x\in \tilde{\Lambda}}P(\varphi_x\geq 0\,\big|\,{\cal F}_{\partial
B(x)})\cdot1_{\Omega_\Lambda^+\cap F}\right]+E\left[\prod_{x\in \tilde{\Lambda}}P(\varphi_x\geq 0\,\big|\,{\cal F}_{\partial B(x)})\cdot1_{\Omega_\Lambda^+\cap F^c}\right],$$
where $F=E_{\delta,\kappa}\cup F_\delta.$ The following lemma shows that we can estimate $\prod_{x\in \tilde{\Lambda}}P(\varphi_x\geq 0\,|\,{\cal F}_{\partial B(x)})$ uniformly
on $F:$

\begin{lemma}\label{estimates}
Let $0<\gamma<1.$ The following hold:
\begin{itemize}
\item[(a)] For $L$ large enough, there exist a constant $c_1>0$ such that \begin{eqnarray}\label{Edelta} E\left[\prod_{x\in \tilde{\Lambda}}P(\varphi_x\geq 0\,\big|\,{\cal
F}_{\partial B(x)})\cdot 1_{\Omega_\Lambda^+\cap
E_{\delta,\kappa}}\right]\leq
\exp\left(-c_1 N^{d-2k+\gamma}\right).\end{eqnarray}
\item[(b)] For $N$ large enough, there exists a constant $c_2$ such that \begin{eqnarray}\label{Fdelta}E\left[\prod_{x\in \tilde{\Lambda}}P(\varphi_x\geq 0\,\big|\,{\cal
F}_{\partial B(x)})\cdot 1_{\Omega_\Lambda^+\cap
F_{ \delta}}\right]\leq \exp\left(-c_2N^{d-2k+\gamma}\right).\end{eqnarray}
\end{itemize}
Both constants depend on $L, \theta$ and $\delta$ but not on $N$.
\end{lemma}

\begin{Proof}
In both cases, we use standard estimates on the centred Gaussian variables $m_x-\varphi_x$ under $P(\,\cdot\,|{\cal
F}_{\partial B(x)}).$ \\
\it(a) \rm Since $G_L\longrightarrow G,$ we have that, for $L$ large enough, $4k(G-\kappa)/2G_L\leq 2k-\gamma.$ We therefore get on $E_{\delta,\kappa}$  
\begin{eqnarray*}
\prod_{x\in \tilde{\Lambda}}P(\varphi_x\geq 0\,\big|\,{\cal F}_{\partial B(x)})&\leq& P\left(\varphi_0-m_0\leq a_N\,\big|\,{\cal F}_{\partial
B(x)}\right)^{\delta B}\\
&\leq &\left(1-\frac{\sqrt{G_L}}{a_N}\exp\left(-\frac{a_N^2}{2G_L}\right)\right)^{\delta B}\\
&\leq&\exp\left(-c_1N^{d-2k+\gamma}\right).
\end{eqnarray*}
\it(b) \rm On $F_{\delta}$ we have for some constants $C>0, c_2>0,$ and for $N$ large enough
\begin{eqnarray*}
\prod_{x\in \tilde{\Lambda}}P(\varphi_x\geq 0\,\big|\,{\cal F}_{\partial B(x)})&\leq& \sum_{\lambda=0}^\infty P(\varphi_0-m_0\geq
N^{\lambda}
\,\big|\,{\cal F}_{\partial B(x)})^{\delta N^{-2k -(2-\gamma)\lambda +\gamma}B}\\
&\leq&\sum_{\lambda=0}^\infty \left(\exp\left(-\frac{N^{2\lambda}}{2G_L}\right)\right)^{\delta N^{-2k -(2-\gamma)\lambda +\gamma}B}\\
&\leq&\sum_{\lambda=0}^\infty\exp\left(-CN^{d-2k+\gamma}\right)^{N^{\gamma\lambda}}\\
&\leq&\exp\left(-c_2 N^{d-2k+\gamma}\right).
\end{eqnarray*}
\end{Proof}

Thus we only need to consider $F^c,$ where we can easily bound  $\sum_{x\in\tilde{B}_r}m_x.$ Write
$$\sum_{x\in\tilde{B}_r}m_x=\sum_{x:m_x>a_N}m_x+\sum_{x: -1<m_x\leq
a_N}m_x+\sum_{\lambda=0}^{\lambda_{\max}}\sum_{x:-N^{\lambda+1} <m_x\leq
N^{\lambda}}m_x$$ and bound the three parts separately: On $E_{\delta,\kappa}^c,$ at least $(1-\delta)$ of the $m_x$ are at height at least
$a_N,$ so for the first part we get
\begin{eqnarray}\label{p1}\sum_{m_x>a_N}m_x\geq(1-\delta)B\,a_N\,.\end{eqnarray} The second term can be estimated easily by writing
\begin{eqnarray}\label{p2}\sum_{-1<m_x\leq a_N}m_x\geq -B\,.\end{eqnarray} 
Finally, since on $F_\delta^c$ there is
$$\big|\{x\in\tilde{B}_r:-N^{\lambda+1}<m_x\leq -N^{\lambda}\}\big|\leq \big|\{x\in\tilde{B}_r:m_x\leq -N^{\lambda}\}\big|\leq\frac{\delta}{N^{2k+(2-\gamma)\lambda-\gamma}}B,$$
we get

\begin{eqnarray}\label{p3}\sum_{\lambda=0}^{\lambda_{\max}}\left[\sum_{ -N^{\lambda+1}<m_x\leq N^{\lambda}}m_x\right]&\geq&
-\sum_{\lambda=0}^{\lambda_{\max}}B\cdot \delta \cdot N^{-2k-(2-\gamma)\lambda+\gamma}\cdot N^{\lambda+1}\nonumber\\
&=&-B\cdot \delta \cdot N^{-2k+\gamma+1} \sum_{\lambda=0}^{\lambda_{\max}} N^{-(1-\gamma)\lambda}\nonumber\\
&\geq& -c\cdot B\cdot N^{-2k+\gamma+1}.
\end{eqnarray}

The three estimates (\ref{p1}), (\ref{p2}), (\ref{p3}) together give

\begin{eqnarray}\label{m_x}\frac{1}{B}\sum_{x\in\tilde{B}_r}m_x\geq (1-\delta) a_N+O(1)\,\end{eqnarray}
on $F^c.$ Let $f_r\geq 0\, \,(r\in I).$ Then (\ref{m_x}) implies
\begin{eqnarray}
P(\Omega_N^+\cap F^c)&\leq& P\left(\sum_{r\in I}f_r\frac{1}{B}\sum_{x\in \tilde{B}_r}m_x>(1-\delta)a_N\sum_{r\in
I}f_r+O(1)\right)\nonumber\\
&\leq&\exp\left(\frac{-(1-\delta)^2 a_N^2(\sum_{r\in I}f_r)^2+O(\sqrt{\log N})}{2 \text{var}(\sum_{r\in I}f_r\frac{1}{B}\sum_{x\in
\tilde{B}_r}m_x)}\right).\label{genest}
\end{eqnarray}
Now we can conclude the proof of the upper bound as in \cite{BDZ}. Since $m_x$ is the conditional expectation
$E(\varphi_x\,|\,{\cal F}_{\partial B(x)})=E(\varphi_x\,|\,{\cal F}_{\Lambda}),$ we have by Jensen's inequality
$$\text{var}\left(\sum_{r\in I}f_r\frac{1}{B}\sum_{x\in \tilde{B}_r}m_x\right)\leq \text{var}\left(\sum_{r\in I}f_r\frac{1}{B}\sum_{x\in
\tilde{B}_r}\varphi_x\right).$$
Define $f_\theta:\R^d\to \R$ by $f_\theta(t)=\sum_{r\in I}f_r1_{A_r}(t).$ One easily sees that
$$\sum_{r\in I}f_r\frac{1}{B}\sum_{x\in
\tilde{B_r}}\varphi_x=\frac{1}{B}\sum_{x\in \tilde{\Lambda}}f_\theta\left(\frac{x}{N}\right)\varphi_x \ \text{ and }\quad \sum_{r\in I}f_r=\frac{1}{B}\sum_{x\in
\tilde{\Lambda}}f_\theta\left(\frac{x}{N}\right)\,,$$ and consequently $$\text{var}\left(\frac{1}{B}\sum_{r\in I}f_r\sum_{x\in
\tilde{B_r}}\varphi_x\right)=\frac{1}{B^2} \sum_{x,y\in \tilde{\Lambda}}f_\theta\left(\frac{x}{N}\right)f_\theta\left(\frac{y}{N}\right)G(x,y)\,.$$

Thus we obtain, using Lemma \ref{estimates} and (\ref{genest}),
\begin{equation*}\begin{split}
\limsup_{N\to \infty}&\frac{1}{N^{d-2k}\log N}\log P(\Omega_N^+)\\
&\leq\limsup_{N\to \infty}\frac{1}{N^{d-2k}\log
N}\log P(\Omega_\Lambda^+\cap F^c)\\
&\leq\limsup_{N\to \infty}\frac{1}{N^{d-2k}\log N}\frac{-(1-\delta)^2 4k(G-\kappa)\log N(\sum_{r\in
I}f_r)^2}{2\text{var}(\sum_{r\in I}f_r\frac{1}{B}\sum_{x\in \tilde{B_r}}m_x)}\\
&=-(1-\delta)^2 2k(G-\kappa)\limsup_{N\to\infty}\frac{1}{N^{d-2k}}\cdot \frac{(\sum_{x\in
\tilde{\Lambda}}f_\theta(\frac{x}{N}))^2}{\sum_{x,y \in
\tilde{\Lambda}}f_\theta(\frac{x}{N})f_\theta(\frac{y}{N})G(x,y)}\,.\\
\end{split}\end{equation*}

As in \cite{BDZ}, the proof is now concluded by applying Proposition \ref{konstanten}, taking the supremum over all possible $f_\theta$ and
letting $\kappa \to 0$ and $\delta \to 0.$\\

\section{Proof of the height estimate}\label{section_height}

To prove Theorem \ref{thmheight}, there are two directions to show. The first was proved in Theorem 2.2 of \cite{Sakagawa}: For any $\varepsilon>0, \eta>0,$ and $ z\in V_N,$ such
that $V_{N,\varepsilon}(z)\subset V_N,$ 

\begin{eqnarray}
\lim_{N\to\infty}P\left(\frac{\overline{\varphi}_{N,\varepsilon}(z)}{\sqrt{\log
N}}\leq\sqrt{4kG}-\eta\;\Bigg|\;\Omega_N^+\right)&=&0\,.\label{unten}
\end{eqnarray}

We will now use Theorem \ref{thmer} to show the other bound:

\begin{proposition}\label{oben}
For any $\varepsilon>0, \eta>0$ and $z\in V_N,$ with $V_{N,\varepsilon}(z)\subset V_N$
\begin{eqnarray}
\lim_{N\to\infty}P\left(\frac{\overline{\varphi}_{N,\varepsilon}(z)}{\sqrt{\log
N}}\geq\sqrt{4kG}+\eta\;\Bigg|\;\Omega_N^+\right)&=&0\,.
\end{eqnarray}
\end{proposition}

The proof for the lattice free field in \cite{BDZ} uses the FKG-inequality for the conditional measure, which does not hold in our case.
Similarly to Section 2, we can handle this problem by carefully estimating the probability that, on $\Omega_N^+,$ the local sample mean
of the field is higher than $\sqrt{4kG\cdot \log N}.$ This is done by comparing
\newcommand{\phimittel}{\overline{\varphi}_{N,\varepsilon}}$\phimittel (z)$ with the average of the conditional expectations $m_x.$

\begin{proof}
First, let $z=0,$ set $\phimittel:=\phimittel (z),$ and $V_{N,\varepsilon}:=V_{N,\varepsilon}(0).$ Fix $L$ as in Section 2 and recall
the definition of the subgrid $\overline{\Lambda},$ the boxes
$B(x)$ and their $K-$boundary $\partial B(x).$ In this section, $\tilde{\Lambda}$ denotes the set $\{x\in \overline{\Lambda}:\partial
B(x)\subset V_{N,\varepsilon}\},$ and $\Lambda=\cup_{x\in\tilde{\Lambda}}\partial B(x).$ For $r\in \R^d$ and $0<\theta<1$ let $A_r$ be
defined as in section 2, and set $I=\{r\in \theta\mathbbm{Z}^d:A_r\subset V_\varepsilon\}, \partial A_r\cap \partial V_\varepsilon =
\emptyset\},$ where $V_\varepsilon=[-\varepsilon,\varepsilon]^d.$  Set $B_r=NA_r,$ and $\tilde{B}_r=B_r\cap\tilde {\Lambda}.$ As before,
set $m_x:=E(\varphi_x\,|\,{\cal{F}}_{\partial B(x)}) $ for $x\in \tilde{B}_r.$ 

We want to estimate 
$$P\left(\phimittel\geq(\sqrt{4kG}+\eta)\sqrt{\log{N}}\,\big|\,\Omega_N^+\right)=\frac{1}{P(\Omega_N^+)}\,P\left(\{\phimittel\geq(\sqrt{4kG}+\eta)\sqrt{\log{N}}\}\cap\Omega_N^+\right)\,.$$
Recall $F$ from Section 3, fix $t>0$ and set $D_t:=\{\varphi:\text{there is} \, r\in I\, \text{such that} \,\frac{1}{|\tilde{B}_r|}\sum_{x\in\tilde{B}_r}(\varphi_x-m_x)<-t\}.$ Then we can write
\begin{equation*}\begin{split}
P\left(\{\phimittel\geq(\sqrt{4kG}+\eta)\sqrt{\log{N}}\}\cap\Omega_N^+\right)=\,&\,P\left(\{\phimittel\geq(\sqrt{4kG}+\eta)\sqrt{\log{N}}\}\cap\Omega_N^+\cap F\right)\\
+\,&\,P\left(\{\phimittel\geq(\sqrt{4kG}+\eta)\sqrt{\log{N}}\}\cap\Omega_N^+\cap F^c\cap D_t\right)\\
+\,&\,P\left(\{\phimittel\geq(\sqrt{4kG}+\eta)\sqrt{\log{N}}\}\cap\Omega_N^+\cap F^c\cap D_t^c\right)\,.
\end{split}\end{equation*}

We have seen in the last section that the first term is negligible compared to $P(\Omega_N^+).$ For the second part, recall that conditioned on ${\cal{F}}_{\partial B(x)}$,
the $\varphi_x-m_x,\, x\in \tilde{B}_r,$ are independent centred Gaussian variables with variance $G_L.$ Thus for the variance of the average we get
$${\rm var}\left(\frac{1}{|\tilde{B}_r|}\sum_{x\in\tilde{B}_r}(\varphi_x-m_x)\;\Bigg|\,{\cal F}_\Lambda \right)=\frac{1}{|\tilde{B}_r|^2}\sum_{x\in\tilde{B}_r}{\rm
var}\left(\varphi_x-m_x\,\big|\,{\cal{F}}_{\partial B(x)}\right)=\frac{1}{|\tilde{B}_r|}\cdot G_L\,.$$ 
We can therefore find constants $c_1>0,$ and $c_2=c_2(\theta)>0$ such that
\begin{eqnarray}\label{links}P\left(D_t\cap F^c\cap \Omega_N^+\right)&\leq&c_2 E\left(P\left(\frac{1}{|\tilde{B}_r|}\sum_{x\in\tilde{B}_r}(\varphi_x-m_x)<-t\;\Bigg|\,{\cal F}_\Lambda\right)\cdot
1_{F^c\cap\Omega_N^+}\right)\nonumber\\
&\leq&c_2 \exp\left(\frac{-t^2\cdot c_1N^d}{2G_L}\right),
\end{eqnarray}
which is also negligible compared with $P(\Omega_N^+).$ Therefore we only need to estimate
$$\limsup_{N\to\infty}\frac{1}{P(\Omega_N^+)} P\left(\{\phimittel\geq(\sqrt{4kG}+\eta)\sqrt{\log{N}}\}\cap
\Omega_N^+\cap F^c\cap D_t^c\right)\,.$$

For this purpose we bound $\frac{1}{|\tilde{B}_r|}\sum_{x\in \tilde{B}_r}\varphi_x$ from below on $\Omega_N^+\cap F^c\cap D_t^c.$ Write
\begin{eqnarray}\label{summe}\frac{1}{|\tilde{B}_r|}\sum_{x\in \tilde{B}_r}\varphi_x=\frac{1}{|\tilde{B}_r|}\sum_{x\in
\tilde{B}_r}(\varphi_x-m_x)+\frac{1}{|\tilde{B}_r|}\sum_{x\in \tilde{B}_r}m_x \end{eqnarray}
and recall from the last section, that on $\Omega_N^+\cap F^c$ 
\begin{eqnarray}\label{rechts}\frac{1}{|\tilde{B}_r|}\sum_{x\in \tilde{B}_r}m_x\geq (1-\delta)a_N+O(1)\end{eqnarray}
 for any $\delta >0$ and $\kappa>0.$ This implies that on $\Omega_N^+\cap F^c\cap D_t^c,$ we have
$$\frac{1}{|\tilde{B}_r|}\sum_{x\in \tilde{B}_r}\varphi_x\geq(1-\delta)a_N+ O(1)\,.$$

Since we can repeat this argument with any shift of the subgrid $\Lambda,$ and average over all
shifts, we conclude that on $\Omega_N^+\cap F^c\cap D_t^c$ 
\begin{equation}\label{abschphi}\frac{1}{|B_r|}\sum_{x\in B_r}\varphi_x\geq (1-\delta)a_N+O(1)\,.\end{equation}

\newcommand{\phir}{\overline{\varphi}_r}From now on we will abbreviate $\frac{1}{|B_r|}\sum_{x\in B_r}\varphi_x$ by $\phir.$
For $\kappa'>0,$ set $C_{\kappa'}:=\{\varphi:\;{\text {there exists}}\; r_0 \;{\text {such that}}\;\overline{\varphi}_{r_0} \geq
(\sqrt{4k(G-\kappa})+\kappa')\sqrt{\log N}\}.$ It follows from (\ref{abschphi}) that, on $\Omega_N^+\cap F^c\cap D_t^c,$ for every $\eta>0$ we can find $\kappa'>0$ and $r_0\in
I$ such that, for $N \to \infty,$
\begin{equation*}P\left(\phimittel\geq(\sqrt{4kG}+\eta)\sqrt{\log N}\right)
\leq P\left(\{\phir\geq (1-\delta)a_N\; \forall r\in I\}\cap C_{\kappa'} \right)\,.\end{equation*}
 Let $f_r>0, r\in I.$
\begin{equation*}\begin{split}
P\bigg(\phir\geq&(1-\delta)a_N\quad \forall r\in I, \;  
\overline{\varphi}_{r_0} \geq (\sqrt{4k(G-\kappa})+\kappa')\sqrt{\log N}\bigg)\\
&\leq P\left(\sum_{r\in I}f_r\phir\geq(1-\delta)a_N\cdot\sum_{r\in I}f_r + \kappa' f_{r_0}\sqrt{\log N}\right)\\
&\leq\exp\left(\frac{-\left((1-\delta)a_N\sum_{r\in I}f_r+ \kappa' f_{r_0}\sqrt{\log
N}\right)^2}{2\text{var}\left(\sum_{r \in I}f_r\phir\right)}\right)\,.
\end{split}\end{equation*}
Defining $f_\theta$ as in the last section, we have
$$\sum_{r\in I}f_r=\frac{1}{|B_r|}\sum_{x\in V_{N,\varepsilon}}f_\theta\left(\frac{x}{N}\right)$$
and
\begin{eqnarray*}
\text{var}\left(\sum_{r\in I}f_r\phir\right)&=&\frac{1}{|B_r|^2}\sum_{r\in I}\sum_{s\in I}f_rf_s\sum_{x\in B_r}\sum_{y\in B_s}G(x,y)\\
&=&\frac{1}{|B_r|^2}\sum_{x,y\in V_{N,\varepsilon}}f_\theta\left(\frac{x}{N}\right)f_\theta\left(\frac{x}{N}\right)G(x,y)\\
&=&O(N^{-d+2k})\,.
\end{eqnarray*}

Similarly to the end of Section 2, we can then optimise over $f_\theta,$ use Proposition \ref{konstanten}, and let $\kappa$ and $\delta$ tend to 0. Then we see that
there is a constant $c>0,$ such that
$$\limsup_{N\to\infty}\frac{1}{ N^{d-2k}\log{N}}\log{P\left(\phimittel \geq (\sqrt{4kG}+\eta)\sqrt{\log N}\right)}
\leq-4k\cdot G\cdot C_k(V)-c\,.$$
Now we apply Theorem \ref{thmer}, and obtain
\begin{equation*}
\limsup_{N\to\infty}\frac{1}{ N^{d-2k}\log{N}}\log{P\left(\phimittel\geq (\sqrt{4kG}+\eta)\sqrt{\log N}\;\big|\;\Omega_N^+\right)}\leq-c\,,
\end{equation*}which proves the claim in the case $z=0.$
The case of an arbitrary $z$ is obtained by repeating the same arguments on a shifted grid.
\end{proof}

Theorem \ref{thmheight} now follows immediately from (\ref{unten}) and Proposition \ref{oben}. This proves the height estimate.

\section{Green's function and $k-$Capacity}

In this section, we prove that there are several equivalent expressions for the capacity $C_k(V),$ as in the case $k=1.$  A crucial step is the
decay of the Green's function $G(x,y)$ as $|x-y|$ tends to infinity. In the
case of the free field, the local central limit theorem for the simple random walk yields a decay  of order $|x|^{-d+2}$ (see \cite{BD},
\cite{Lawler}). In our model, we do not have a random walk representation, since the entries of the ``transition matrix'' $I-J$ can be negative. Nevertheless, using the methods
of Section 3 of \cite{LeGall}, one can, without the use of a local central limit theorem, obtain a decay of $G$
of order $|x|^{2k-d}.$ This was done by Sakagawa:

\begin{lemma}\label{lemma1}(\cite{Sakagawa}, Lemma 5.1)
Let $d\geq 2k+1.$ Then there is a constant $\eta_k$ such that
$$\lim_{|x|\to\infty}\frac{G(0,x)}{|x|^{2k-d}}=\frac{1}{q_k}\eta_k\,.$$
\end{lemma}

Define now a function $g_k:\R^d\to\R$ by $g_k(x)=\frac{\eta_k}{q_k}|x|^{2k-d}$ and a positive compact operator $K_k$ on $L^2(V)$ by
$$K_kf(x)=\int_Vg_k(x-y)f(y)dy \quad (x\in V)\,.$$ From the above lemma we get, for $|x-y|\to\infty,$
\begin{eqnarray}\label{greens}|g_k(x-y)-G(x,y)|=o(|x-y|^{-d+2k})\,.\end{eqnarray}

In this section we use the short notation $<f,g>_V:=\int_Vf(x)g(x)dx,$ for suitable functions $f, g$ and $V\subset \R^d.$ Note first the following (see
also \cite{Sakagawa}, Lemma 5.2): 

\begin{lemma} \label{lemma2}Let $d\geq 2k+1.$ For all $h, f\in H^k(V),$
$$\frac{q_k}{(2d)^k}<h, K_k\frac{1}{(2d)^k}(-\Delta)^kf>_V=<h,f>_V.$$
\end{lemma}

\begin{proof}
In order to distinguish between the discrete and the continuous
Laplacian, we denote them by $\Delta_d$ and $\Delta_c$ respectively. Using (\ref{greens}) we obtain
\begin{equation*}\begin{split}
\frac{q_k}{(2d)^k}<h,&K_k(-\Delta_c)^kf>_V\\
&=\lim_{N\to\infty}\frac{1}{N^{2d}}\sum_{x\in V_N}\sum_{y\in
V_N}h\left(\frac{x}{N}\right)\cdot g_k\left(\frac{x}{N}-\frac{y}{N}\right)\cdot \frac{q_k}{(2d)^k}\cdot \left((-\Delta_c)^kf\right)\left(\frac{y}{N}\right)\\
&=\lim_{N\to\infty}\frac{1}{N^{2d}}\sum_{x\in V_N}h\left(\frac{x}{N}\right) \sum_{y\in
V_N}N^{d-2k}\cdot G\left(\frac{x}{N},\frac{y}{N}\right)\cdot \frac{1}{N^{-2k}} q_k \left((-\Delta_d)^kf\right)\left(\frac{y}{N}\right)\\
&=\lim_{N\to\infty}\frac{1}{N^d}\sum_{x\in V_N}h\left(\frac{x}{N}\right)\cdot f\left(\frac{x}{N}\right)
=\frac{1}{q_k}<h,f>_V.
\end{split}\end{equation*}
\end{proof}

We can now prove the equivalence of several expressions for the $k-$th order capacity $C_k(V).$ Proposition \ref{konstanten} below was used implicitly in Section 5 of
\cite{Sakagawa} (Lemma 5.2). As we are not aware of a reference, we include the proof here.

\begin{proposition}\label{konstanten}
Let $V=[-1,1]^d, d\geq 2k+1.$ Then
\begin{equation*}\begin{split}
\inf&\left\{\frac{q_k}{(2d)^k}\int_{\R^d}|(-\nabla)^k h|^2dx: H \in H^k_0(\R^d), h\geq1_V\right\}\\
&=\sup\left\{2<f, 1_V>_V-<f, K_kf>_V: f\in L_2(V)\right\}\\
&=\sup\left\{\frac{<f, 1_V>^2_V}{<f,K_kf>_V}: f\in L_2(V)\right\}.
\end{split}\end{equation*}
\end{proposition}

\begin{proof}
Let us prove the first equality. Notice that $M:=\{h\in H^k_0(\R^d): h\geq1_V\}$ is a closed convex subset of the Hilbert space $H_0^k(\R^d),$ and
thus has a minimizer $h_0$ for the Sobolev-norm on $H_0^k(\R^d).$ But this means exactly that $h_0$ minimizes $\int_{\R^d}|(-\nabla)^k h|^2dx$ for
$h\in M.$ It is immediate that $h_0=1$ on $V.$ Furthermore, $(-\Delta)^k h_0=0$ outside $V.$ To see this, set $g(\varepsilon)=\int_{\R^d}|(-\nabla)^k
h+\varepsilon \varphi|^2dx $ for some $\varphi\in C_c^\infty(\R^d\setminus V).$ Then $\frac{dg}{d\varepsilon}\Big|_{\varepsilon=0}=0,$
because $h_0$ is a minimizer of the integral. But this implies $<(-\Delta )^kh_0,\varphi>_{\R^d\setminus V}=<(-\nabla)^kh_0,(-\nabla)^k\varphi>_{\R^d\setminus V}=0$ for all $\varphi \in
C_c^\infty(\R^d\setminus V)$ and thus $(-\Delta)^k h_0=0$ on $\R^d\setminus V.$\\
There exist $\tau_n\in C_0^\infty(\R^d), n\in \N$ such that $\lim_{n\to\infty}<h_0-\tau_n,
(-\Delta)^k(h_0-\tau_n)>_{\stackrel{\circ}{V}}=0$ and $\tau_n=h_0$ on $\R^d\setminus V,$ where $\stackrel{\circ}{V}$ denotes the interior of $V.$ Set $$f_n=\frac{q_k}{(2d)^k}(-\Delta)^k\tau_n\,.$$

For every $n,$ $f_n$ belongs to $L_2(\R^d),$ and, by the fact that $f_n=0$ outside $V,$ Lemma \ref{lemma2} and partial integration yield

\begin{eqnarray*}
2<f_n,\tau_n>_V-<f_n, K_kf_n>_V=\frac{q_k}{(2d)^k}<(-\Delta)^k\tau_n,\tau_n>_{\R^d}
=\frac{q_k}{(2d)^k}\int_{\R^d}|(-\nabla)^k\tau_n|^2dx\,.
\end{eqnarray*}

Moreover, as in \cite{BD}, $\lim_{n\to\infty}|<f_n, 1_V-\tau_n>_{L^2(\R^d)}|=0.$ Together with the above this yields
\begin{equation*}\begin{split}
\sup\{2<f,1_V>_V-<f,K_kf>_V\}
&\geq\limsup_{n\to\infty}\{2<f_n,1_V>_V-<f_n,K_kf_n>_V\}\\
&=\frac{q_k}{(2d)^k}\int_{\R^d}|(-\nabla)^kh_0|\,dx\,,
\end{split}\end{equation*}
which gives one direction in the first equation. The other direction is an elementary calculation based on Lemma \ref{lemma2}.\\
The second equation follows by expanding $f$ in a basis of eigenvectors of the compact positive operator $K_k.$ Maximising shows
that both sides are equal to $\sum_{i\in \N}\frac{<e_i,1_V>^2}{\lambda_i},$ where the $e_i$ are the eigenvectors of $K_k$ and
$\lambda_i$ the corresponding eigenvalues.

\end{proof}

\noindent\bf{Acknowledgements.} \rm I wish to express my thanks to Erwin Bolthausen for his help and encouragement. Thanks also to Hironobu Sakagawa for the discussion of a problem in an earlier version.

\bibliographystyle{plain}

\end{document}